\documentstyle[12pt]{article}

\newcommand\qed{\hfill $\Box$ \vspace{5mm}}
\newcommand\pthree{{\rm \bf P\/}^3}
\newcommand\pfour{{\rm \bf P\/}^4}

\newcommand\lm{\lambda}
\newcommand{\binom}[2]
{ \left(\begin{array}{c}#1\\#2\end{array}\right)}

\begin{document}

\title{An improved bound for the degree of smooth surfaces in $\pfour$ not of
general type}
\author{Michele Cook}
\maketitle

This is an addendum to the paper of Braun and Fl{\o}ystad ([BF])
on the bound for the degree of smooth surfaces in $\pfour$ not of general
type. (In fact, any statements made here without comment will
be found there.)  Using their construction and the regularity of curves
in $\pthree$, one may lower the bound a little more.
We will prove the following:

\vspace{5mm}
{\bf Theorem.}

{\it Let $S$ be a smooth surface of degree $d$ in $\pfour$ not of general type.
Then either $d \leq 70$ or $S$ lies on a hypersurface in $\pfour$
of degree 5 and $d \leq 80$. Thus $d \leq 80$. }

\vspace{2mm}
(We will also indicate how one may actually get the bound
down to 76.)

\vspace{3mm}
The main result of [BF] is to bound $\chi {\cal O}_S$ from below using
generic initial ideal theory and
information coming from a generic hyperplane section of $S$.

\vspace{3mm}
Let $C$ be a generic hyperplane section of $S$, then $C$ has
asociated to it
two pieces of information relevant to this situation;
{\it connected invariants}
$\lm_0 > \lm_1 > \dots >  \lm_{s-1} > 0$ and {\it sporadic zeros} which
are defined as follows:

\vspace{5mm}
{\bf Definition}.
Let ${\rm \bf C}[x_0, x_1, x_2, x_3]$ be the ring of polynomials of $\pthree$,
with reverse lexicographical ordering.
Let $C$ be a curve in $\pthree$, then the generic initial ideal
of $C$,  ${\rm gin}(I_C)$,
is generated by elements of the form $x_0^{i}x_1^{j}x_2^{k}$.
we say that a monomial $x_0^ax_1^bx_2^c$ is a {\it sporadic zero}
if $x_0^ax_1^bx_2^c \notin {\rm gin}(I_C)$, but there exists
$c' > c$ such that $x_0^ax_1^bx_2^{c'} \in {\rm gin}(I_C)$.

\vspace{3mm}
Let $\alpha_t$ is the number of sporadic zeros of degree t
and $m$ be the maximal degree of the sporadic zeros.

\vspace{3mm}
If $\pi$ is the genus of $C$, then
$$ \pi =1+ \sum_{i=0}^{s-1}(\binom{\lm_i}{2}+(i-1)\lm_i) -
\sum_{t=0}^{m}\alpha_t$$
and
$$1+ \sum_{i=0}^{s-1}(\binom{\lm_i}{2}+(i-1)\lm_i) \leq
\frac{d^2}{2s}+(s-4)\frac{d}{2}+1-\frac{r(s-r)(s-1)}{2s} = G(d,s),$$

for $d > (s-1)^2+1$ and $d \equiv r$ (mod $s$) $0 \leq r < s.$

(This is due to the work
of Gruson and Peskine on the numerical invariants of points (see [GP]).)

\vspace{5mm}
In [BF] they show that
$$ \chi {\cal O}_S \geq
\sum_{t=0}^{s-1} (\binom{\lm_t+t-1}{3} - \binom{t-1}{3})
- \sum_{t =0}^{m} \alpha_t (t-1), $$
and if $s \geq 2$ and $d > (s-1)^2+1$,
$$\sum_{t=0}^{s-1} (\binom{\lm_t+t-1}{3} - \binom{t-1}{3}) \geq
s\binom{\frac{d}{s}+\frac{s-3}{2}}{3}+1-\binom{s-1}{4}.$$

\vspace{5mm}
They Combine this with
the double point formula (see [H, pg 434])

$$ d^2-5d-10(\pi - 1) + 2(6 \chi - K^2) = 0,$$

and the fact that if the degree of $S$ is $d > 5$ then $K^2 \leq 9$,
to get

$$\begin{array}{lll}
 18 \geq 2K^2 & = & d^2-5d-10(\pi - 1) + 12 \chi  \\
   & \geq & d^2-5d
-10(\frac{d^2}{2s}+(s-4)\frac{d}{2} - \sum_{t=0}^{m}\alpha_t) \\
&& +12( s\binom{\frac{d}{s}+\frac{s-3}{2}}{3}+1-\binom{s-1}{4}
- \sum_{t =0}^{m} \alpha_t (t-1))  \\
   & = & d^2-5d-10(\frac{d^2}{2s}+(s-4)\frac{d}{2}) \\
&& +12( s\binom{\frac{d}{s}+\frac{s-3}{2}}{3}+1-\binom{s-1}{4})
 - \sum_{t =0}^{m} \alpha_t (12t-22). \ \ (\star)
\end{array}$$

\vspace{3mm}
We will use regularity conditions to find an upper bound for

$$ A = \sum_{t=0}^{m} \alpha_t(12t-22).$$

\vspace{3mm}
First we should note that if $s = \min\{k | h^0({\cal I}_S(k) \neq 0\}$.
Then by
the work of Ellingsrud and Peskine ([EP]) if $S$
is a surface, of degree $d$, not of general type then either
$=d \leq 90$ or $s \leq 5$. Furthermore, if one imitates
their argument with $s=6$, one gets either $d \leq 70$ or $s \leq 6$.
We will assume $s \leq 6$. We also know that if $s \leq 3$
then $d \leq 8$. Hence we only need to bound the degree of surfaces
not of general type with $s = 4, 5,\  {\rm or}\  6$.

\vspace{3mm}
By [BS] the regularity of an ideal I is equal to the
regularity of ${\rm gin}(I)$, which is the highest degree of
a minimal generator of ${\rm gin}(I)$.
Furthermore, by [GLP], the regularity of a smooth curve of degree $d$
in $\pthree$ is $\leq d-1$.
Hence the largest degree of the minimal generators of
${\rm gin}(I_C)$ is $\leq d-1$ and so all the sporadic zeros of $C$
are in degree $\leq  d-2$.

\vspace{3mm}
We will now bound the number of sporadic zeros.
Let $\gamma = G(d,s)-\pi$.
Any bound on $\gamma$ will also bound the number of sporadic zeros.
By [EP], $\gamma \leq \frac{d(s-1)^2}{2s}$. Furthermore, by the
double point formula and the fact that if $S$ is a surface not of
general type $K^2 < 6 \chi$, we get $\pi \geq \frac{d^2-5d+10}{10}$ and thus
$\gamma \leq  \frac{d^2}{2s}+(s-4)\frac{d}{2}+1-\frac{d^2-5d+10}{10}$.
Taking the minimum of these bounds for $\gamma$, we get, for
$s=4, \gamma \leq \frac{9d}{8}$, for $s=5, \gamma \leq d$ and
for $s=6, \gamma \leq \frac{d(90-d)}{60} $.

Furthermore, by connectedness of the invariants
$\lambda_0 \leq \frac{d}{s} + s - 1$ and $\lambda_1 \leq \frac{d}{s} + s - 2$.

\vspace{3mm}
Putting all this together, we have

$$
A \leq  \sum_{t=\frac{d}{s}+s-1}^{d-2}(12t-22)
+ \sum_{t=\frac{d}{s}+s-1}^{\gamma - d + \frac{2d}{s}+2s-2}(12t-22)
$$

\vspace{5mm}
and using the bounds on $\gamma$, we get
$$\begin{array}{ll}
{\rm for} \ s=4, & A \leq \frac{243}{32}d^2-9d+192 \\
{\rm for} \ s=5, & A \leq \frac{162}{25}d^2 -16d+300. \\
\end{array}$$

\vspace{5mm}
Substituting back into the original equation $(\star)$ above, we get

$$\begin{array}{lll}
{\rm for} \ s=4, & 0 \geq \frac{1}{8}d^3-\frac{275}{32}d^2+\frac{7}{2}d-195 &
{\rm and \  hence} \ d \leq 68 \\
{\rm for} \ s=5, & 0 \geq \frac{2}{25}d^3-\frac{162}{25}d^2+4d-318 &
{\rm and \  hence} \ d \leq 80.
\end{array}$$

\vspace{3mm}
Now suppose $S$ does not lie on a hypersurface of degree 5. Then
the degree of $S  = d \leq 90$. By [EP], we know that either
$d \leq 70$ or $S$ is contained in a hypersurface of degree 6.
Therefore it remains to show that there are no smooth surfaces not of
general type with $s=6$ and $71 \leq d \leq 90$.

Suppose such a surface existed. Then $\gamma \leq \frac{71(90-71)}{60}$,
i.e. $\gamma \leq 22$ and so
$$ A \leq \sum_{t = \frac{d}{6}+5}^{\frac{d}{6}+5+21}(12t-22)
= 44d+3608.$$

Substituting as above we get

$$ 0 \geq \frac{d^3}{18}+\frac{2}{3}d^2-\frac{119}{2}d-\frac{7321}{2}.$$

But this is a contradiction for $71 \leq d \leq 90$. \qed

\vspace{3mm}
By more careful consideration of the equations involved, it is
posssible to lower the bound a little more:

Making use of the equation
$$ 0 \leq \sum_{t=0}^{m}\alpha_t \leq
\sum_{i=0}^{s-1}(\binom{\lm_i}{2}+i(\lm_i)) - \frac{d^2+5d}{10}$$
arising from the double point formula, and
$$\begin{array}{ll}
 0 \geq & d^2+5d-18-10(\sum_{i=0}^{s-1}(\binom{\lm_i}{2}+i(\lm_i))\\
& +12(\sum_{i=0}^{s-1}(\binom{\lm_i+i-1}{3})-\binom{s-1}{4})
 - \sum_{t =0}^{m} \alpha_t (12t-22)
\end{array}$$
it is possible to show that $d \leq 76$ and $s \leq 8$.
This is done be writing out all the possible connected invariants
of length $s$ and degree $d$
and checking the inequalities on ``Mathematica''. (I would like to thank
Rich Liebling for sending me the program which calculates the invariants.)

\vspace{3mm}
The temptation at this point is to try to lower the bound
by considering each degree individually and using ad hoc methods. For
example, we know for certain configurations of the invariants the curve
is Arithmetically Cohen-Macaulay and hence has no sporadic zeros.
But, to make any significant improvements using these methods,
one needs to understand which configurations of the sporadic zeros can
occur.

\vspace{3mm}
{\it acknowledgements}. I would like to thank Mei-Chu Chang for asking
me to talk on the paper of Braun and Fl{\o}ystad, which made me read
the details of the paper. I would like to thank Robert Braun
and Gunnar Fl{\o}ystad for answering all my questions and Robert for
pointing out the better bound for $\gamma$ for $s=6$.  Finally, I would like
to thank Bruce Crauder and Sheldon Katz for making me feel welcome
during my stay here at OSU.

\vspace{5mm}

{\bf References}

{\bf [BS]} D. Bayer, M. Stillman {\it A criterion for
detecting m-regularity}, Invent. Math. 87 (1987) pp1-11.

{\bf [BF]} R. Braun, G. Fl{\o}ystad {\it A bound for the degree
of smooth surfaces in $\pfour$ not of general type}, Composito
Mathematica, Vol. 93, No. 2, September(I) 1994 pp211-229.

{\bf [EP]} G. Ellingsrud, C. Peskine {\it Sur les surfaces lisses de
$\pfour$.}, Invent. Math. 95 (1989) pp1-11.

{\bf [GLP]} L. Gruson, R. Lazarsfeld, C. Peskine {\it On a Theorem of
Castelnuovo, and the equations defining space curves}, Invent. Math.
72 (1983) pp491-506.

{\bf [GP]} L. Gruson, C. Peskine {\it Genres des courbes de l'espace
projectif}, Lecture Notes in Mathematics, Algebraic Geometry,
Troms{\o} 1977, 687 (1977), pp31-59.

{\bf [H]} R. Hartshorne {\it Algebraic Geometry}, Springer-Verlag (1977).

\vspace{5mm}
Michele Cook

Department of Mathematics

Oklahoma State University

Stillwater, OK 74078.

e-mail mcook\verb+@+math.okstate.edu

\end{document}